\documentclass[12pt]{article}

\usepackage{amsmath}
\usepackage{amssymb} 
\usepackage{latexsym}

\usepackage{txfonts}

\newcommand{\mcal}[1]{{\mathcal {#1}}} %%CAPITALS only!

\newtheorem{theorem}{Theorem}  %% [section]
\newtheorem{lemma}[theorem]{Lemma}
\newtheorem{corollary}[theorem]{Corollary}
\newtheorem{proposition}[theorem]{Proposition}

{\end{enumerate}}

\newenvironment{proof}[1]{\smallskip \noindent {\bf #1}}{\qed\smallskip}

\def\qed{\ifhmode\unskip\nobreak\fi\ifmmode\ifinner\else\hskip5pt\fi\fi
 \hfill\hbox{\hskip5pt\vrule width4pt height6pt depth1.5pt\hskip1pt}}

\newcommand{\ov}[1]{\ensuremath{\overline{#1}}}

\newcommand{\Int}{\ensuremath{\operatorname{\rm Int}}}

\newcommand{\sgn}{\ensuremath{\operatorname{\rm{sign}}}}

\newcommand{\dimply}{\ensuremath{\:\Longleftrightarrow\:}}

\newcommand{\R}[1]{\ensuremath{{\mathbb R}^{#1}}} %bb superscrip R

\newcommand{\Rp}{\ensuremath{{\mathbb R}_+}}   %bb R_plus
\newcommand{\ZZ}{\ensuremath{{\mathbb Z}}}	   %bold Z
\newcommand{\Np}{\ensuremath{{\mathbb N}_{+}}}  %bold N^plus
\newcommand{\om}[1]{\ensuremath{\omega(#1)}}

\newcommand{\pde}[2]
	{\ensuremath{\frac{\partial #1}{\partial #2}}}	%%partial derivative 
\newcommand{\ode}[2]{\ensuremath{\frac{d#1}{d#2}}}    	%% derivative

\def\d#1dt{\frac{d#1}{dt}}    %%variable ODE left side
\newcommand{\lam}{\ensuremath{\lambda}}
\newcommand{\Lam}{\ensuremath{\Lambda}}
\newcommand{\del}{\ensuremath{\delta}}

\newcommand{\sig}{\ensuremath{\sigma}}
\newcommand{\Sig}{\ensuremath{\Sigma}}
\newcommand{\Gam}{\ensuremath{\Gamma}}
\newcommand{\gam}{\ensuremath{\gamma}}

\newcommand{\co}{\colon\thinspace} %% Colon with correct spacing for
\newcommand{\sk}{\smallskip\noindent}
     
\newcommand{\skk}{\medskip\noindent}

\newcommand{\p}{\ensuremath{\partial}}

\reversemarginpar       %%PUTS LABELS IN LEFT MARGIN
 \def\mylabel#1{\label{#1}}
\title{\mbox{Attractors in coherent systems of differential equations}}

\author{David Angeli, Dip. di Sistemi e Informatica, University of Firenze\\
Morris W. Hirsch, Dep. of Mathematics, University of California Berkeley\\
Eduardo D. Sontag, Dep. of Mathematics, Rutgers University}

\def\E{{\mathcal E}}

\def\Ko{{\mathsf K^o}}

\def\J{{\mathsf J}}
\def\Jo{\mathsf J^o}

\newcommand{\comment}[1]{}

\begin{document}

\maketitle
% \begin{flushright} {\footnotesize /Cascades/d5MH.tex \normalsize}
% \end{flushright}

\begin{abstract} 

Attractors of cooperative dynamical systems are particularly simple; for
example, a nontrivial periodic orbit cannot be an attractor.
This paper provides characterizations of attractors for the wider class of
\emph{coherent} systems, defined by the property that no directed feedback
loops are negative.
Several new results for cooperative systems are obtained in the process.

\end{abstract}
 \tableofcontents

%%%%%%%%%%%%%%%%%%%%%%%%%%%%%%%%%%%%%%%%%%%%%%%%%%%%%%%%%%%%
\section{Introduction}   \mylabel{sec:intro}
%%%%%%%%%%%%%%%%%%%%%%%%%%%%%%%%%%%%%%%%%%%%%%%%%%%%%%%%%%%%

We consider differential equations 
%%%%%%%%%%%%%%%%%%%%%%%%%%%%%%%%%%%%%%%%%%%%%%%%%%%%%%%%%%%%
\begin{equation}		\label{eq:basic}
%%%%%%%%%%%%%%%%%%%%%%%%%%%%%%%%%%%%%%%%%%%%%%%%%%%%%%%%%%%%
\ode x t = F(x),  \qquad x\in X,\quad t \ge 0,
\end{equation}
where $X\subset \R n$ is convex, its interior is dense in $X$, and the
vector field $F\co X\to \R n$ extends to a $C^1$ vector field  on
an open set.   The maximally defined solutions $t\mapsto \Phi_t
(a), t\ge 0, a\in X$ generate the local semiflow $\Phi:=\{\Phi_t\}_{t\in\Rp}$.
We refer to  $F$ (or $(F, X, \R n)$, or $(F, X,
\R n,\Phi)$) as a {\em system}.  Dynamical notions
are applied interchangeably to $F$ and $\Phi$.

Many biological situations are modeled by {\em cooperative} systems:
$\pde{F_j}{x_i}\ge 0$ if $j\ne i$.  The biological interpretation is
that an increase of species $i$ tends to increase the population
growth rate of every other species $j$.  In this case $\Phi$ is {\em
monotone}, meaning it preserves the vector ordering.
This causes the crude dynamics of a cooperative system to be
comparatively simple; for example, there are no attracting cycles and
every orbit is nowhere dense (Hadeler \& Glas \cite {HadelerGlas83},
Hirsch \cite {Hirsch82a}).

Here we show that some of the dynamical advantages of cooperative
systems extend to systems having a significantly weaker property: $F$
is {\em coherent} 
(another name is \emph{positive feedback system})
if
whenever $i_0, \dots, i_\nu,\; \nu \in \{1,\dots,n\}$ are such that
\[
   i_\nu=i_0,\quad i_{k-1}\ne  i_k \textrm{ and } \frac{\partial F_{i_{k-1}}}{\partial x_{i_k} }  \nequiv 0 \quad (1 \leq k \leq \nu) \]
then,
\[   \pde{F_{i_{k-1}}}{x_{i_k}}(x) \text{ does not change sign}\qquad 
   (1\le k \le \nu)
\]   
and
%%%%%%%%%%%%%%%%%%%%%%%%%%%%%%%%%%%%%%%%%%%%%%%%%%%%%%%%%%%%
\begin{equation}		\label{eq:coherent}
%%%%%%%%%%%%%%%%%%%%%%%%%%%%%%%%%%%%%%%%%%%%%%%%%%%%%%%%%%%%
\pde{F_{i_0}}{x_{i_1}} (x)\ \cdots \
 \pde{F_{i_{\nu-1}}}{x_{i_\nu}}(x) \ge  0,\qquad (\forall \, x\in X). 
\end{equation}

%
%When formulas for the $F_i$ are given, the validity of cooperativity
%and coherence depends on the geometry of $X$ in $\R n$.  For example,
%suppose $X\subset \R 2$ and $F\co X\to \R 2$ is given by
%\[
%\begin{split}
%F_1(x_1, x_2) &= -x_1 x_2 + g(x_2)\\
%F_2 (x_1,x_2) &=  -x_1^2  + h(x_1) 
%\end{split}
%\]
%with $C^1$ functions $g$ and $h$.  This system is cooperative when $X$
%lies in the left halfplane, and coherent when $X$ lies in the first or
%third quadrant, but not otherwise.
 
Our chief combinatorial result, Theorem \ref{th:cascade}, shows that
by permuting the variables $x_i$ and changing the signs of some of
them, any coherent system can be transformed into a dynamically
equivalent system $(F,X, \R n, \Phi)$ with the following properties: 
\begin{itemize}
\item $F$ is not merely coherent, it  has the stronger property of being {\em
  quasicooperative}:   for any  $(i_1, \dots,i_m)$ as above,
  each factor in the left hand side of (\ref{eq:coherent}) is $\ge
  0$

\item if $F$ is not cooperative, there exists a cooperative system
 $(F^1,X^1, \R {n_1}, \Phi^1)$, $1\le n_1< n$, such that the the
 natural projection 
%%%%%%%%%%%%%%%%%%%%%%%%%%%%%%%%%%%%%%%%%%%%%%%%%%%%%%%%%%%%
%\begin{equation}		\label{eq:Pidef}
%%%%%%%%%%%%%%%%%%%%%%%%%%%%%%%%%%%%%%%%%%%%%%%%%%%%%%%%%%%%
\[
  \Pi \co \R n \to \R {n_1}, \quad (x_1,\ldots,x_n) \mapsto
 (x_1,\ldots,x_{n_1})
\]
%\end{equation}
%%
maps $X$ onto $X^1$ and semiconjugates $F$ to $F^1$ and $\Phi$ to $\Phi^1$:
%%
%%%%%%%%%%%%%%%%%%%%%%%%%%%%%%%%%%%%%%%%%%%%%%%%%%%%%%%%%%%%
%\begin{equation}		\label{eq:cascade}
%%%%%%%%%%%%%%%%%%%%%%%%%%%%%%%%%%%%%%%%%%%%%%%%%%%%%%%%%%%%
\[
  \Pi\circ F (x) = F^1\circ \Pi (x), \qquad \Pi\circ \Phi_t (x) =
  \Phi_t^1\circ \Pi (x) \ \; \text {if $\Phi_t (x)$ is defined}
\]
%% end{equation}
\end{itemize}
Mild geometrical conditions on $X$ guarantee that for each equilibrium
$p$ of $F^1$, the restriction of $F$ to $X_p:=X\cap\Pi^{-1}(p)$ is
equivalent to a quasicooperative  system $(\hat F_p, \hat X_p, \R {n-n_1})$.
This is the basis for inductive proofs of our main results. 

\smallskip
We turn to our main topic, attractors.  An {\em attractor} for $F$ 
is a nonempty invariant continuum $A\subset X$ that uniformly
attracts all points in some neighborhood of $A$.  If the attraction is
not necessarily uniform we talk instead of an {\em attracting set}.
Three types of attractors $A$ have received special attention:

\smallskip
{\em Point attractors}: $A$ is single point, necessarily an equilibrium. 

\smallskip {\em Periodic attractors}:  $A$ is a cycle, i.e., a periodic orbit
that is not an equilibrium.

\smallskip {\em Strange attractors}, often called ``chaotic''.  This
somewhat vague term signifies that $A$ is neither an equilibrium nor a
cycle, and usually that $A$ is topologically transitive and exhibits
``sensitive dependence on initial conditions''.  Some authors also
require that periodic orbits be dense in $A$. 

This paper is motivated by the question: What kind of nonequilibrium
attractors $A$ can exist in coherent systems?   Theorem \ref{th:A} shows that 
$A$ cannot be topologically transitive; Theorems \ref{th:B} and
\ref{th:C} give further dynamical information.  Other  results apply
to more general monotone local semiflows.  

%%%%%%%%%%%%%%%%%%%%%%%%%%%%%%%%%%%%%%%%%%%%%%%%%%%%%%%%%%%%
\subsection{Statement of results}   \mylabel{sec:oultine}
%%%%%%%%%%%%%%%%%%%%%%%%%%%%%%%%%%%%%%%%%%%%%%%%%%%%%%%%%%%%

A set is {\em finitely transitive} for a system (or a local semiflow)
if it is the union of the omega limit sets of finitely many of its
points. 

%%%%%%%%%%%%%%%%%%%%%%%%%%%%%%%%%%%%%%%%%%%%%%%%%%%%%%%%%%%%
\begin{theorem}		\mylabel{th:A}
%%%%%%%%%%%%%%%%%%%%%%%%%%%%%%%%%%%%%%%%%%%%%%%%%%%%%%%%%%%%
 A finitely transitive attracting set $A$ for a  system $(F,
 X, \R n)$ reduces to an equilibrium in the following cases:
\begin{description}

\item[(i)]  $F$ is coherent, and $X$ is open in
  $\R n$ or relatively open in a coordinate half-space

\item[(ii)] $F$ is quasicooperative, and every point of $A$ is
   strongly accessible in $X$ from above, or every point of $A$ is
   strongly accessible in $X$ from below

\item[(iii)]  $F$ is cooperative, and each point of $A$ is
   strongly accessible in $X$ from above or below

\end{description}
\end{theorem}
%%%%%%%%%%%%%%%%%%%%%%%%%%%%%%%%%%%%%%%%%%%%%%%%%%%%%%%%%%%%
A stronger conclusion, Theorem \ref{th:fintrans},  holds for cooperative
systems.  
  
The following  result requires no additional geometrical conditions on
$X$: 
%%%%%%%%%%%%%%%%%%%%%%%%%%%%%%%%%%%%%%%%%%%%%%%%%%%%%%%%%%%%
\begin{theorem}		\mylabel{th:B}
%%%%%%%%%%%%%%%%%%%%%%%%%%%%%%%%%%%%%%%%%%%%%%%%%%%%%%%%%%%%
If $(F, X, \R n)$, $n\ge 2$ is a coherent system, every  orbit is
nowhere dense.
\end{theorem}
%%%%%%%%%%%%%%%%%%%%%%%%%%%%%%%%%%%%%%%%%%%%%%%%%%%%%%%%%%%%
\paragraph{Conjecture}  
{\em In a coherent system with $n>1$, every orbit closure has measure
zero.}  Even for cooperative systems this is known only for $n=2$.

An attractor is {\em global} if it attracts all points of $X$.  An
equilibrium is {\em globally asymptotically stable} if it is the
global attractor.  The following theorem needs $X$ to be open:
%%%%%%%%%%%%%%%%%%%%%%%%%%%%%%%%%%%%%%%%%%%%%%%%%%%%%%%%%%%%
\begin{theorem}		\mylabel{th:C}
%%%%%%%%%%%%%%%%%%%%%%%%%%%%%%%%%%%%%%%%%%%%%%%%%%%%%%%%%%%%
Let $(F, X, \R n)$ be a coherent system with $X$ open in $\R n$.
Assume there exists a global attractor $A$.  Then there exists an
equilibrium, and if it is unique it is globally asymptotically
stable.
\end{theorem}
%%%%%%%%%%%%%%%%%%%%%%%%%%%%%%%%%%%%%%%%%%%%%%%%%%%%%%%%%%%%
Proposition \ref{th:moncon} extends a basic result
previously known only for strongly order-preserving local semiflows. 
The development of the concept ``attractor'' is discussed in the
Appendix.

%%%%%%%%%%%%%%%%%%%%%%%%%%%%%%%%%%%%%%%%%%%%%%%%%%%%%%%%%%%%
\subsection{Motivations}
%%%%%%%%%%%%%%%%%%%%%%%%%%%%%%%%%%%%%%%%%%%%%%%%%%%%%%%%%%%%

A coherent system is one whose interaction graph (defined below)
has no directed negative loops.  A more restrictive condition, for
graphs that are not necessarily strongly connected, is
the requirement that the graph has no \emph{un}directed negative loops:
in that case, one may always perform an elementary change of variables 
(defined below) that transforms such a system into a cooperative one.
In a classical and often-quoted 1981 paper, R. Thomas conjectured that
coherent systems do not have any periodic attractors: ``the presence of at
least one negative loop in the logical structure appears as a necessary (but
not sufficient) condition for a permanent periodic behavior'' \cite{thomas81}.
It has often been claimed (see e.g. \cite{pigolotti2007}) that Thomas'
conjecture was settled in \cite{snoussi,gouze98}.  However, these references
only dealt with the more restricted monotone case.
Theorem \ref{th:A} in this paper settles the question.
%%%%
We refer the reader to \cite{almostmonotone_journal} for further comments on
the relevance of these concepts to molecular systems biology, and to
\cite{pf07} for numerical simulations which suggest that systems
that are ``close'' to having the coherence property might have, in some
statistical sense, simpler attractors.

%%%%%%%%%%%%%%%%%%%%%%%%%%%%%%%%%%%%%%%%%%%%%%%%%%%%%%%%%%%%
\subsection{Structure of proofs}
%%%%%%%%%%%%%%%%%%%%%%%%%%%%%%%%%%%%%%%%%%%%%%%%%%%%%%%%%%%%
The proofs of Theorems \ref{th:A}, \ref{th:B} and \ref{th:C} have a
common pattern which we now discuss.  Let $\mcal T$ stand for one of
these theorems.  It is proved first for a cooperative system, which
includes the case $n=1$.  The proof proceeds by induction on $n$.  A
coherent system which is not cooperative is transformed, by permuting
and changing signs of 
variables, to a system $(F, X, \R n)$ having the following properties:

\begin{itemize}

\item $F$ is quasicooperative

\item there is a system 
 $(F^1, X^1, \R {n_1})$ with $n_1 < n$, such that the
natural projection $\Pi \co \R n \to \R {n_1}$
satisfies 
%%%%%%%%%%%%%%%%%%%%%%%%%%%%%%%%%%%%%%%%%%%%%%%%%%%%%%%%%%%%
\begin{equation}		\label{eq:Pi}
%%%%%%%%%%%%%%%%%%%%%%%%%%%%%%%%%%%%%%%%%%%%%%%%%%%%%%%%%%%%
\Pi(X)=X^1, \qquad  \Pi\circ F (x)  = F\circ \Pi (x), \quad (x\in X)
\end{equation}

\item $F^1$ is cooperative
\end{itemize}
It follows that $\Pi$ semiconjugates  the local semiflow $\Phi$ of $F$
to the local semiflow $\Phi^1$ of $F^1$:
%%%%%%%%%%%%%%%%%%%%%%%%%%%%%%%%%%%%%%%%%%%%%%%%%%%%%%%%%%%%
\begin{equation}		\label{eq:piphi}
%%%%%%%%%%%%%%%%%%%%%%%%%%%%%%%%%%%%%%%%%%%%%%%%%%%%%%%%%%%%
 \Pi\circ \Phi_t(x) = \Phi^1_t\circ \Pi(x) \ \; \mbox{if $\Phi_t(x)$ is defined}  
\end{equation}
We summarize this by saying that  
$\Pi \co (F, X, \R n)\twoheadrightarrow (F^1, X^1, \R {n_1})$ (or
$\Pi \co F \twoheadrightarrow F^1$) is a {\em cascade}.   We also
allow the {\em  trivial} cascade, for which  $F= F^1$. 
\smallskip

For each equilibrium $p$ of $F^1$ the affine subspace
$E_p:=\Pi^{-1} (p)$ is a coset of the kernel of $\Pi$.  The {\em
  canonical chart} 

%%%%%%%%%%%%%%%%%%%%%%%%%%%%%%%%%%%%%%%%%%%%%%%%%%%%%%%%%%%%
\begin{equation}		\label{eq:canonchart}
%%%%%%%%%%%%%%%%%%%%%%%%%%%%%%%%%%%%%%%%%%%%%%%%%%%%%%%%%%%%
T_p\co E_p\approx \R {n-n_1}, \quad 
(x_1,\ldots ,x_n)\mapsto (x_{n_1+1}, \ldots, x_n)
\end{equation}
is an affine automorphism.

The vector field $F$, being tangent to $E_p$ along $X_p$, restricts to
a vector field $F_p$ in $X_p:=X_p:=X\cap E_p$, and $\Phi$ restricts to
a local semiflow $\Phi_p$ in $X^p$.  The hypothesis of $\mcal T$ will
ensure that the relative interior of $X_p$ in $E_p$ is dense in $X_p$.
The canonical chart converts the {\em fibre system} $(F_p, X_p, E_p)$ into
a system $(\hat F_p, \hat X_p,\R {n-n_1})$.

{\em 
We identify each fibre system $(F_p, X_p, E_p)$ with 
$(\hat   F, \hat X_p,\R {n-n_1})$ 
by means of the canonical chart.} 
 Thus $F_p$ has an interaction graph $\Gam (F_p):=\Gam (\hat F_p)$.
  We ascribe to $F_p$ the property of being cooperative,
  quasicooperative or coherent whenever that property holds for $\hat
  F_p$.

Theorem $\mcal T$ holds for the cooperative system $F^1$, and it holds
for all fibre systems by the inductive assumption.  The induction is
completed by showing that this implies $\mcal T$ also holds for $(F,
X, \R n)$.

There is a delicate point regarding the domains of these systems.  The
proofs for cooperative systems use special properties of
$X$, such as every point being strongly accessible from above.
These properties are postulated in the hypotheses of the main
theorems.  To make the induction work, the same properties must be
verified for the systems obtained by elementary coordinate changes,
and also for fibre systems.  This means that the class of domains $X$
referred to in the theorems must
be preserved by permuting and changing signs of variables, and by
intersecting $X$ with the affine subspaces $E_p$.  For this reason $X$
is usually required to be an open set in $\R n$ or a relatively open
subset of a coordinate halfspace.  

%%%%%%%%%%%%%%%%%%%%%%%%%%%%%%%%%%%%%%%%%%%%%%%%%%%%%%%%%%%%
\subsection{Local semiflows}
%%%%%%%%%%%%%%%%%%%%%%%%%%%%%%%%%%%%%%%%%%%%%%%%%%%%%%%%%%%%
 A {\em local semiflow} $\Phi$ in a metrizable space $Z$ is a collection
$\Phi=\{\Phi_t\}_{t\in \Rp}$ of continuous maps  $\Phi_t\co D_t\to R_t$
between nonempty subsets of $Z$, with $D_t$ open. 
The notation $\Phi_t x$ indicates $x\in D_t$, absent
contraindications.  $\Phi$ is required to have the following properties:
\begin{itemize}

\item The set $\Omega:=\{(t,x)\in \Rp\times Z\co x\in D_t\}$ is an open
neighborhood of $\{0\}\times Z$  in $\Rp\times Z$, and the map 
$\Omega \to Z, \quad (t,x)\mapsto \Phi_t x$
is continuous. 

\item  $x\in (\Phi_s)^{-1} D_t \implies \Phi_s\circ\Phi_t (x)=\Phi_{s+t}(x)$

\item $\Phi_0$ is the identity map of $Z$. 
\end{itemize} 
We also  say that $(\Phi, Z)$ is a local semiflow.
When $\Phi$ is obtained by solving Equation (\ref{eq:basic}) each
map $\Phi_t$ is a homeomorphism, but this is not assumed for general
local semiflows.

%%DA the meaning of the above sentence is not clear. Isn't this
%% always the case ?]}

The {\em orbit}  and {\em omega limit set} of $x$ are respectively
\[\gam(x):=\{\Phi_t x\co x\in D_t\}, \qquad
\om x :=\, \bigcap_{t\ge 0} \ov{\gam (\Phi_t x)}
\]  

$p$ is an {\em equilibrium} if $\Phi_tp=p$ for all $t$.  The set of
equilibria is denoted by $\E (\Phi)$, and by $\E (F)$ when $\Phi$ is
generated by the vector field $F$.

%%%%%%%%%%%%%%%%%%%%%%%%%%%%%%%%%%%%%%%%%%%%%%%%%%%%%%%%%%%%
\subsection{Attractors and attracting sets}   \mylabel{sec:attrac}
%%%%%%%%%%%%%%%%%%%%%%%%%%%%%%%%%%%%%%%%%%%%%%%%%%%%%%%%%%%%
 We call $A$ {\em positively
  invariant} for if $\Phi_t (a)$ is 
defined and belongs to $A$ for all $t\ge 0$, $a\in A$, and 
{\em invariant}  if in addition $A$ is nonempty and $\Phi_t (A)=A$ for
all $t\ge 0$.  We say that  $A$ 
{\em attracts} $x$ if  $\ov{\gam(x)}$ is compact
and  $\om x\subset A$.  The set of such points $y$ is the {\em basin}
of $A$.  

$A$ is {\em topologically transitive} if it is the omega limit set of
one of its points, and {\em finitely transitive } if it is the union of
the omega limit sets of finitely many of its points.

We call  $A$  {\em attracting} if it is invariant,
connected and  compact, and its basin is a neighborhood of $A$. 
If in addition $A$ has arbitrarily small positively invariant neighborhoods,
 $A$ is an {\em attractor}.\footnote
{
There are many definitions of ``attractor'' in current use,  not
mutually consistent.  The one adopted 
here is equivalent to that of Conley \cite{Conley78}, and (for compact
invariant sets) those of Hale \cite{Hale88} and Sell \& You \cite{SY}.
It is analogous to the definitions for discrete-time systems in  
Smale \cite{Smale67} and Akin \cite{Ak93}.
}

%%%%%%%%%%%%%%%%%%%%%%%%%%%%%%%%%%%%%%%%%%%%%%%%%%%%%%%%%%%%
\subsection{Ordered spaces}  
%%%%%%%%%%%%%%%%%%%%%%%%%%%%%%%%%%%%%%%%%%%%%%%%%%%%%%%%%%%
By an {\em ordered space}  we mean a topological space
$Z$ together  an order  relation $\mcal R\subset
Z\times Z$ that is topologically closed.  If $x,y \in Z$ we write:
%%%%%%%%%%%%%%%%%%%%%%%%%%%%%%%%%%%%%%%%%%%%%%%%%%%%%%%%%%%%
\begin{equation}		\label{eq:order}
%%%%%%%%%%%%%%%%%%%%%%%%%%%%%%%%%%%%%%%%%%%%%%%%%%%%%%%%%%%%
\mbox{$x\succeq y$ and $y\preceq x$ if $(x,y)\in \mcal R$, \ \ \
$x\succ y$ and $y\prec x$ if $x\succeq y, x\ne y$}
\end{equation}
\smallskip
The {\em
vector order} in any subspace of $\R n$ is  defined by 
\[u \succeq v
  \dimply u-v \in \Rp^n
\]
where $\Rp^n$ denotes the the positive orthant $[0,\infty)^n\subset\R
n$. 

A subset of an ordered space is {\em unordered} if none of its
points are related by $\succ$. 

Every subspace  $X\subset Z$ inherits an order relation from $Z$.  
If $M\subset Z$ then
$x\succ M$ means $x\succ y$ for all $y \in M$,
and similarly for the other relations in (\ref{eq:order}). 
For $x, y \in X$ we write 
\[
\begin{split} 
x \rhd_X  y   &\text{ if } \  x\succ N,\ y\in \Int_X (N),\\
x \lhd_X  y   &\text{ if } \   x\prec N,\ y\in \Int_X (N) 
\end{split}
\]
%%DA: specified N needs to be open.. otherwise it boils down to \succ, just take N={ y }
for some open $N \subset X$.
Note the notational anomaly that $x \rhd_X  y$ and $y \lhd_X x$ are not
equivalent statements for general ordered spaces.  They are
equivalent, however, if $X\subset \R n$ is open and has  the vector
ordering.   For example, in $X=\Rp^2$ we have $(0,0)  \lhd_X  (0,1)$ but
$(0,1) \not  \rhd_X (0,0) $.

%%DA: remark above should be illustrated by means of a simple
%%  example

%%
Let $X$ be a subset of an ordered space $Z$. We call $q\in X$  {\em
strongly accessible in $X$ from above (respectively, from below)}  if
every neighborhood of $q$ in $X$ contains a point $x \rhd_X q$
(respectively, 
$x \lhd_X q$).\footnote
{Slightly stronger properties with the same names are used in Hirsch \&
 Smith \cite {HirschSmith05}.}

All our results are valid when $X$ is an open set in $\R n$, and some
are valid for  special kinds of nonopen sets, especially  open subsets of a
{\em coordinate halfspace}  of $\R n$, which means a set  
\[
   \{x\in \R n\co \alpha x_l \ge  c_l\}
\]
for some choice of $l\in\{1,\dots,n\}$, $\alpha \in\{\pm 1\}$, 
$(c_1,\dots,c_n)\in \R n$. 
We rely on the following fact, whose proof is left to the reader:
%%%%%%%%%%%%%%%%%%%%%%%%%%%%%%%%%%%%%%%%%%%%%%%%%%%%%%%%%%%%
\begin{lemma}		\mylabel{th:access}
%%%%%%%%%%%%%%%%%%%%%%%%%%%%%%%%%%%%%%%%%%%%%%%%%%%%%%%%%%%%
Assume $X\subset \R n$ has the vector ordering.  If $X$ is an open
subset of $\R n$, or a relatively open subset of a coordinate
halfspace, every point of $X$ is strongly accessible from above and
below in $X$. \qed
\end{lemma}
%%%%%%%%%%%%%%%%%%%%%%%%%%%%%%%%%%%%%%%%%%%%%%%%%%%%%%%%%%%%
Note also that if $X$ is an open subset of $\Rp^n$, all points of $X$
are strongly accessible  in $X$ from above.

%%%%%%%%%%%%%%%%%%%%%%%%%%%%%%%%%%%%%%%%%%%%%%%%%%%%%%%%%%%%
\section{Cascades}   \label{sec:casc}
%%%%%%%%%%%%%%%%%%%%%%%%%%%%%%%%%%%%%%%%%%%%%%%%%%%%%%%%%%%%
Let $(F, X, \R n)$ and $(F^1, X^1, \R {n_1})$ be systems with $1\le
n_1 <n$ and assume $\Pi \co F\twoheadrightarrow F^1$ is a cascade (see
(\ref{eq:Pi})).   
This implies 
%%%%%%%%%%%%%%%%%%%%%%%%%%%%%%%%%%%%%%%%%%%%%%%%%%%%%%%%%%%%
\begin{equation}		\label{eq:fi}
%%%%%%%%%%%%%%%%%%%%%%%%%%%%%%%%%%%%%%%%%%%%%%%%%%%%%%%%%%%%
\pde{F_i}{x_j} =0\ \ \text{if } \  i \le n_1 <j,
\qquad (i, j\in\{1,\dots,\nu)
\end{equation}
and the Jacobian matrices of $F$ have lower triangular block
decompositions of the form
%%%%%%%%%%%%%%%%%%%%%%%%%%%%%%%%%%%%%%%%%%%%%%%%%%%%%%%%%%%%
\begin{equation}		\label{eq:block}
%%%%%%%%%%%%%%%%%%%%%%%%%%%%%%%%%%%%%%%%%%%%%%%%%%%%%%%%%%%%
F'(x) =\left [
	\begin{array}{ll}
 	M_{11} (x) 	& O    	     		\\
	M_{21} (x) 	& M_{22}(x) 
	\end{array}
	\right ]
\end{equation} 
where  $M_{11} (x) = (F^1)' (\Pi x)\in \R {n_1\times n_1}$, and $O$
stands for 
a matrix of zeroes.  
The following diagrams commute for each $t\ge 0$:

\hspace*{1.7in}
\begin{picture}(120,100)(5, -5)
\put(10,10){$D(\Phi^1_t)$}		%lower left
\put(110,10){$R(\Phi^1_t)$}		%lower right
\put (9,60){$D(\Phi_t)$}		%upper left
\put(110,60){$R(\Phi_t)$}		%upper right

%% Arrows and labels

\put(50,13){\vector(1,0){51}}		%bottom
	 \put(69, 18){$\Phi^{1}_t$}
\put (50,63){\vector (1,0){50}}		%top
	 \put (69, 68){$\Phi_t$}
\put (23,50) {\vector (0,-1){25}}	%left
	 \put (6, 35){$ \Pi$}
\put(125,50){\vector(0,-1){25}}		%right
	 \put(130,35){$ \Pi$}
\end{picture}

For  $p\in \E (F^1)$ let $T_p\co E_p\approx \R{n-n_1}$ be the canonical
chart.  Set $\hat X_p=T_p (X_p)$ and define $\hat F_p\co \hat X_p \to
\R{n-n_1}$ to be the unique vector field transformed by
$(T_p)^{-1}$ to $F_p$, that is, 

%%%%%%%%%%%%%%%%%%%%%%%%%%%%%%%%%%%%%%%%%%%%%%%%%%%%%%%%%%%%
\begin{equation}		\label{eq:Tp}
%%%%%%%%%%%%%%%%%%%%%%%%%%%%%%%%%%%%%%%%%%%%%%%%%%%%%%%%%%%%
 \hat F (T_p x ) =  T_p F (x), \qquad (x\in X_p)
\end{equation}
The local semiflows of $F_p$ and $\hat F_p$ are conjugate under $T_p$.
For $(\hat F_p, \hat X_p, \R {n-n_1})$ to be a system it is necessary
and sufficient that the relative interior of $X_p$ in $E_p$ be dense
in $E_p$.  When this holds we call  $(F_p, X_p, E_p)$ the  {\em
fibre system} over $p$ and identify it with 
$(\hat F_p, \hat X_p, \R {n-n_1})$ by $T_p$, 

The interaction graph $\Gam(F_p):=\Gam (\hat F_p)$ is determined by
the signs of the entries in the block $M_{22}(x)$ in (\ref{eq:block}).
The next lemma gives convenient conditions ensuring this.

  Consider the
following conditions:
\begin{itemize}
\item C$_1 (X, \R n)$:\ $X$ is open in $\R n$
\item C$_2 (X, \R n)$:\  $X$ is open in a coordinate halfspace of $\R n$
\item C$_3 (X, \R n)$:\  $X$ is open in $\Rp^n$
\item C$_4 (X, \R n)$:\  $X$ is a rectangle
\end{itemize}
%%%%%%%%%%%%%%%%%%%%%%%%%%%%%%%%%%%%%%%%%%%%%%%%%%%%%%%%%%%%
\begin{lemma}		\mylabel{th:fibsys}
%%%%%%%%%%%%%%%%%%%%%%%%%%%%%%%%%%%%%%%%%%%%%%%%%%%%%%%%%%%%
Assume  $\Pi \co F\twoheadrightarrow F^1$ a cascade as above and $p\in
\E (F^1)$.   Suppose $p\in \E (F^1)$, and  {\em C$_d (X, \R n)$},
 is satisfied for some $d\in \{1,2,3,4\}$.  Then $(\hat F_p, \hat X_p, \R
 {n-n_1})$ is a system, and  
{\em C$_d (\hat X_p, \R {n-n_1})$} holds.
\end{lemma}
%%%%%%%%%%%%%%%%%%%%%%%%%%%%%%%%%%%%%%%%%%%%%%%%%%%%%%%%%%%%
\begin{proof}  
The verification that  C$_d (X, \R n)$ implies C$_d (\hat X_p, \R
{n-n_1})$, and also that relative interior of $X_p$ in $E_p$ is dense in
$X_p$, is straightforward. 
\end{proof}

%%%%%%%%%%%%%%%%%%%%%%%%%%%%%%%%%%%%%%%%%%%%%%%%%%%%%%%%%%%%
\section{Graphs}   \mylabel{sec:graphs}
%%%%%%%%%%%%%%%%%%%%%%%%%%%%%%%%%%%%%%%%%%%%%%%%%%%%%%%%%%%%
By a {\em directed graph} $\Gam:=(V_\Gam, E_\Gam)$ we mean a nonempty
finite set $V:=V_\Gam$ (the set of vertices) together with a binary relation
$E:=E_\Gam\subset V\times V$ (the set of {\em directed edges}, usually
referred to simply as ``edges'').  We always assume $E$ is totally
nonreflexive i.e., $(i,i)\notin E$.   

 An {\em isomorphism} between a pair of  directed
graphs is a bijection $f$ between their vertex sets such that $f\times
f$ restricts to a bijection $f_*$ between their edge sets.   

Our chief tool for analyzing the crude dynamics of systems $(F, X, \R
n)$ is the {\em interaction graph} $\Gamma:=\Gamma(F)$.  This is the
labeled directed graph with vertex set is $V=V(\Gam):=\{1,\dots,n\}$, whose
set of (directed) edges is
\[
 E=E(\Gam):=\{(j,i)\in V\times V\co j\ne i \ \text{ and }
 \  \frac{\p F_i}{\p x_j} \text{ is
 not identically} \ 0\ \text{in } X \}
\]
Edge $(j,i)$ is assigned the label
$h(j,i)\in\{+1, -1, \theta\}$	according to the rule:
%%%%%%%%%%%%%%%%%%%%%%%%%%%%%%%%%%%%%%%%%%%%%%%%%%%%%%%%%%%%
\begin{equation}		\label{eq:h}
%%%%%%%%%%%%%%%%%%%%%%%%%%%%%%%%%%%%%%%%%%%%%%%%%%%%%%%%%%%%
h(j,i) =  \begin{cases}
  \ \ 1 & \text{if \ $\frac{\p F^i}{\p x^j}(x)\ge 0$\ \ for all } x \in X,\\   
  -1  & \text{if  \ $\frac{\p F^i}{\p x^j} (x) \le 0$\ \ for all } x \in X,\\
       \ \ \theta  & \text{otherwise}
	\end{cases}
\end{equation}
and is respectively called positive, negative or ambiguous.  
A  loop is {\em positive} if each of
its edges is labeled $+1$ or $-1$ and the product of these labels is 
$+1$. 
 
We define three types of graphs in increasing order of generality:

\smallskip
\hspace*{5mm} $\Gam$ is {\em positive} if every edge is positive,

\smallskip
\hspace*{5mm} $\Gam$ is {\em
  quasipositive} if every loop has only positive edges,

  \smallskip
\hspace*{5mm} $\Gam$  has the {\em positive loop property} 
%%
%% and  is called {\em   loop-positive} 
%%
 if every loop is positive.

\sk Paraphrasing some of the earlier definitions, we define corresponding
types of systems $F$ in terms of $\Gam (F)$:

\smallskip
\hspace*{5mm} $F$ is {\em cooperative} if $\Gam (F)$ is positive

\smallskip
\hspace*{5mm} $F$ is {\em quasicooperative} if $\Gam (F)$ is quasipositive, 

\smallskip
\hspace*{5mm} $F$ is {\em coherent}  if $\Gam$ has the positive
  loop property  

\skk Evidently  {\em cooperative $\implies$ quasicooperative $\implies$ 
coherent.} 

\smallskip

The term ``graph'' is shorthand for ``finite directed graph having
edges labeled in $\{1, -1,\theta\}$.''  Graphs are denoted by Greek
capitals $\Gam, \Lam$, perhaps with indices.  The sets of vertices and
edges of $\Gam$ are denoted by $V (\Gam)$ and $E(\Gam)$, respectively,
and the labeling function is denoted by $h_\Gam\co V (\Gam)\to \{1,
-1,\theta\}$.  Two graphs $\Gam, \Lam$ are {\em isomorphic} if there
there is an isomorphism $f\co V(\Gam)\to V (\Lam)$ between the
underlying directed graphs such that $h_\Lam
\circ f_* = h_\Gam$.   
%DA: ^^^^^^^^^^^^^  fixed little typo in the definition of labeling isomorphism

$\Lam$ is a {\em subgraph} of  $\Gam$ provided
 \[
   V(\Lam)\subset V(\Gam), \quad
   E(\Lam)\subset E(\Gam),  \quad h_{\Lam}=h_\Gam|E(\Lam),
 \]
We abuse notation and denote this by  $\Lam
 \subset\Gam$, saying that $\Lam$ contained in $\Gam$.  

If $\Lam,\Lam'$ are
 subgraphs their {\em graph union} is the subgraph with vertex set 
 $V(\Lam)\cup V(\Lam')$ and edge set $E(\Lam)\cup E (\Lam')$. 

A {\em path}  of length $k\in \Np$ is a sequence $(u_0,
\dots,u_k)$ of vertices  such that $(v_{j-1}, v_j)$ is an edge for
  $j=1,\dots, k$. 
The  {\em concatenation}  of an ordered pair $(\lam,\mu)$ of paths,
\[
 \lam =(u_0, \dots,u_k),\quad \mu=(u_k,\dots, u_{k+l}),
\]
is the path 
\[ \lambda\cdot \mu:=(u_0, \dots,u_k, u_{k+1},\dots, u_{k+l})
\] 
obtained by transversing first $\lambda $ and then
$\mu $.  
 
A {\em loop} of length $\mu\in\Np$ is a
sequence of $\mu\ge 2$ edges having the form 
\[ 
  (i_0, i_1), (i_1, i_2), \dots, (i_{\mu-1}, i_\mu), \quad i_\mu=i_0
\]
As our graphs are totally nonreflexive, there are no self-loops:
 $i_j\ne i_{j-1}, \ j=1,\dots, \mu$.  

A loop is {\em positive} (respectively, {\em negative}) if each of
its edges is labeled $1$ or $-1$ and the product of these labels is 
$+1$ (respectively, $-1$).  All other loops are {\em ambiguous}.  

In the next three definitions the labeling plays no role.
A graph is called:

\smallskip
\hspace*{.7cm}  {\em connected } if for each pair of distinct vertices $j,
k$ there is a sequence of vertices $j=i_0, \dots, i_m =k$, $m\in \Np$
such that $(i_{l-1}, i_l)$ or $(i_l, i_{l-1})$ is an edge of $\Lam$,
$(l=1,\dots, m)$  

\smallskip
\hspace*{.7cm} {\em strongly connected } if for any ordered pair
$(a,b)$ of distinct vertices there is a path in $\Lam$ from $a$ to
$b$,

\smallskip
\hspace*{.7cm}  {\em primary } if every  edge belongs to a loop,   

\skk  
These definitions imply:
\begin{itemize}
\item A graph with no edges is primary, but  a graph with only
  one edge is not primary.  

\item The graph union of primary subgraphs  is primary.     

\item
A strongly connected subgraph is primary, and a primary connected
subgraph having more than one  vertex is strongly connected. 
If $\Gam$ is quasipositive, every primary subgraph is positive.
\end{itemize}

A subgraph $\Lam\subset \Gam$ is called:
\begin{itemize}
\item 
{\em full } provided  it contains all edges in
 $\Gam$ joining vertices of $\Lam$,

\item {\em initial } if no edge of $\Gam$ is directed from a vertex outside
$\Lam$ to a vertex of $\Lam$,  

\item{\em terminal } if if no directed edge of $\Gam$
  joins a vertex of $\Lam$ to  a vertex not in  $\Lam$,

\item  {\em fundamental }  if 
is  connected, primary and initial, and no other subgraph
containing $\Lam$ has these properties.  
\end{itemize}

%%%%%%%%%%%%%%%%%%%%%%%%%%%%%%%%%%%%%%%%%%%%%%%%%%%%%%%%%%%%
\begin{lemma}		\mylabel{th:subgraphs}
%%%%%%%%%%%%%%%%%%%%%%%%%%%%%%%%%%%%%%%%%%%%%%%%%%%%%%%%%%%%
The following hold for all subgraphs:

\begin{description}

\item[(a)] fundamental subgraphs are full

\item[(b)] if fundamental subgraphs $\Lam_1, \Lam_2$ share a vertex,
they coincide 

\item[(c)] every connected, primary, initial subgraph is
contained in a unique fundamental subgraph

\end{description}
\end{lemma}
%%%%%%%%%%%%%%%%%%%%%%%%%%%%%%%%%%%%%%%%%%%%%%%%%%%%%%%%%%%%
\begin{proof} (a) and (b) follow directly from definitions. 
(c) is proved by showing that the graph
union of a maximal nested family of connected, primary, initial
subgraphs is fundamental.
\end{proof}

%%%%%%%%%%%%%%%%%%%%%%%%%%%%%%%%%%%%%%%%%%%%%%%%%%%%%%%%%%%%
\subsection{Graphs and systems}   \mylabel{sec:graphsys}
%%%%%%%%%%%%%%%%%%%%%%%%%%%%%%%%%%%%%%%%%%%%%%%%%%%%%%%%%%%%
Let $(F,X,\R n)$ be a system. 

%%%%%%%%%%%%%%%%%%%%%%%%%%%%%%%%%%%%%%%%%%%%%%%%%%%%%%%%%%%%
\begin{proposition}		\mylabel{th:graphs}
%%%%%%%%%%%%%%%%%%%%%%%%%%%%%%%%%%%%%%%%%%%%%%%%%%%%%%%%%%%%
If $\Pi \co F\twoheadrightarrow F^1$ is  a cascade having a fibre
system $F_p$,   then:
\begin{description}

\item[(a)]  $\Gam (F^1)$ is a  full subgraph
of $\Gam (F)$.

\item[(b)] $\Gam (F_p)$ is isomorphic to a  subgraph of
$\Gam (F)$. 
 
\item[(c)] when $F$ is cooperative, quasicooperative or
coherent,  $F^1$ and $F_p$ have the same property.

\end{description}
\end{proposition}
%%%%%%%%%%%%%%%%%%%%%%%%%%%%%%%%%%%%%%%%%%%%%%%%%%%%%%%%%%%%
\begin{proof}  (a) and (b), which imply (c),  are proved by inspecting
  the block   decomposition (\ref{eq:block}) of 
  the matrix of functions $F'(x)$.
\end{proof}

%%DA: I think one assumption is missing, namely that $\Gamma^1$ is
%%  full. I added details accordingly
%%%%%%%%%%%%%%%%%%%%%%%%%%%%%%%%%%%%%%%%%%%%%%%%%%%%%%%%%%%%
\begin{proposition}		\mylabel{th:fund0}
%%%%%%%%%%%%%%%%%%%%%%%%%%%%%%%%%%%%%%%%%%%%%%%%%%%%%%%%%%%%
Let $\Gam^1\subset\Gam (F)$ be an initial full subgraph such that
$V (\Gam_1)=\{1,\dots,n_1\}$.  Then:
\begin{description}

\item[(i)]   there is a cascade 
$  \Pi\co (F, X,\R n) \twoheadrightarrow (F^1,X^1,\R {n_1})$
such that $\Gam (F^1)=\Gam^1$.
 
\item[(ii)]  
When $F$ is quasicooperative,  $F^1$ and all fibre systems
are quasicooperative, and if $\Gam_1$ is primary then $F^1$ is cooperative

\end{description}
\end{proposition}
%%%%%%%%%%%%%%%%%%%%%%%%%%%%%%%%%%%%%%%%%%%%%%%%%%%%%%%%%%%%
\begin{proof} Initiality and fullness of $\Gam^1$ means that 
 (\ref{eq:fi}) holds.  Therefore (\ref{eq:Pi}) defines a cascade
 satisfying (i).  The first assertion in (ii) follows from Proposition
 \ref{th:graphs}(c).  The second assertion holds because $\Gam_1$ is
 quasipositive, and if it is primary  all its edges are in
 loops and hence are positive.
 \end{proof}

%%%%%%%%%%%%%%%%%%%%%%%%%%%%%%%%%%%%%%%%%%%%%%%%%%%%%%%%%%%%
\subsection{Spin assignments}   \mylabel{sec:spin}
%%%%%%%%%%%%%%%%%%%%%%%%%%%%%%%%%%%%%%%%%%%%%%%%%%%%%%%%%%%%

A {\em spin assignment} for a graph $\Gam$ is any function $\sig\co V(\Gam)\to
\{\pm1\}$.  It is {\em consistent} if $h(u,v)=\sig(u)\sig(v)$ for
every edge $(u,v)$ belonging to a loop.
(This terminology is not the same as in \cite{almostmonotone_journal}, where
it was required that every edge be consistent.  With that stronger requirement,
the theorem given below would become a characterization of monotonicity with
respect to an orthant order, a more restrictive property than coherence.)

%%%%%%%%%%%%%%%%%%%%%%%%%%%%%%%%%%%%%%%%%%%%%%%%%%%%%%%%%%%%
\begin{theorem} \mylabel{th:spin} 
%%%%%%%%%%%%%%%%%%%%%%%%%%%%%%%%%%%%%%%%%%%%%%%%%%%%%%%%%%%%
$\Gam$ has  the positive loop property if and only if it has a
consistent spin assignment.
\end{theorem}
%%%%%%%%%%%%%%%%%%%%%%%%%%%%%%%%%%%%%%%%%%%%%%%%%%%%%%%%%%%%
\begin{proof}
Assume $\Gam$ has the positive loop property.  
Let $\Gam'$ be
obtained from $\Gam$ by keeping the same vertices  but deleting the
edges not contained in loops. 
Clearly $\Gam'$ has the positive loop property, and if $\sig$
is a consistent spin assignment on $\Gam'$ it is also consistent on
$\Gam$.  Therefore we can assume every edge $e$  belongs to a
loop and is thus positive. 

{\em Claim: } If $\lam^1,  \lam^2$ are paths from $a$ to $b$ then
$h(\lam^1) = h(\lam^2) \in \{\pm 1\}$. To see this, choose a path $\mu$
from $b$ to 
$a$, which can be done because each edge belongs to a
loop.   Since every loop is positive by hypothesis, for $j=1,2$ we have
\[
  1=h(\lam^j\cdot\mu)=h(\lam^j)h(\mu)
\]
Therefore $h(\lam^1)=h(\mu)=h(\lam^2)$. 

Now fix a vertex $p$ of $\Gam$ and for each vertex $v$ choose a path
$\lam_v$ from $p$ to $v$.  Define $\sig(p)=1$ and $\sig(v)=h(\lam_v)$,
which by the claim is independent of the choice of $\lam_v$.  For any
edge $e=(u,v)$ we can fix $\lam_u$ and define $\lam_v:=\lam_u\cdot
e$. Then have:
\[
 \sig(u)=h(\lam_u),\qquad \sig (v)=h(\lam_u\cdot  e) =h(\lam_u)h(e),
\]
which implies $h(e)=\sig(u)\sig(v)$.   
The converse implication is left to the reader. 
\end{proof}

%%%%%%%%%%%%%%%%%%%%%%%%%%%%%%%%%%%%%%%%%%%%%%%%%%%%%%%%%%%%
\paragraph{Remark}
%%%%%%%%%%%%%%%%%%%%%%%%%%%%%%%%%%%%%%%%%%%%%%%%%%%%%%%%%%%%
The foregoing proof can be expressed homologically.  Let $\hat \Lam$
denote the 1-dimensional cell complex corresponding to a prime
subgraph $\Lam\subset \Gam$ having the vertices of $\Lam$ for
$0$-cells and the directed edges of $\Lam$ for $1$-cells.  In the
cellular chain groups of $\hat \Lam$ with coefficients in $\ZZ_2$
(identified with the multiplicative group $\{\pm 1\}$), a labeling $h$
is a $1$-cochain, spin assignments are $0$-cocycles, and a spin
assignment $\sig$ is consistent for $h$ if its coboundary is $\del
\sig =h$.
As the evaluation of cochains on chains induces a dual pairing $H^1
(\hat \Lam; \ZZ_2)\times H_1 (\hat \Lam; \ZZ_2) \to \ZZ_2$, the
positive loop property  makes the cohomology class of $h$ trivial.
Thus $h=\del\sig$, proving that  $\sig$ is consistent.

\smallskip
A change of variables $x\mapsto y$ is called {\em elementary} if there is
a permutation $i\mapsto i'$ of $\{1,\dots,n\}$ and an $n$-tuple
$\rho\in \{\pm 1\}^n$ such that  $y_i = \rho_i x_{i'}$.  

%%%%%%%%%%%%%%%%%%%%%%%%%%%%%%%%%%%%%%%%%%%%%%%%%%%%%%%%%%%%
\begin{theorem}		\mylabel{th:cascade}
%%%%%%%%%%%%%%%%%%%%%%%%%%%%%%%%%%%%%%%%%%%%%%%%%%%%%%%%%%%%
If a system is coherent,  there is an elementary change of variables
transforming it to a quasicooperative system admitting a cascade over
a cooperative system for which all fibre systems are quasicooperative.
\end{theorem} 
%%%%%%%%%%%%%%%%%%%%%%%%%%%%%%%%%%%%%%%%%%%%%%%%%%%%%%%%%%%%
\begin{proof}
Assume $(F, X,\R n)$ is a coherent system, which by Theorem
\ref{th:spin} has a consistent spin assignment $\sig$.   The
elementary change of variables  $L\co\R n\to\R n$,
\[
  y=Lx, \quad y_i:=\sig(i)x_i
\]
transforms $(F, X,\R n)$ into a system 
\[
  (G, L(X), \R n), \quad L\circ G ***
\] 
such that $\Gam (G)$ and $\Gam (F)$ have the same undirected edges.
For every directed edge $(j,i)$ of $\Gam (G)$:
\[\textstyle
  h_{\Gam (G)}(j,i)=\sgn (\pde{G_j} {y_i})= \sig_j \sig_i \sgn
   (\pde{F_j} {x_i}) \ \sig_j \sig_i h_{\Gam (F)}(j,i)
\]
 If $(j,i)$ belongs to a loop then
$h_{\Gam (F)}(j,i) =\sig_j \sig_i$ by the consistency condition.  Therefore
\[\textstyle
   \sgn (\pde{G_j} {y_i}) (\sig_j \sig_i)^2 = (\pm 1)^2 = 1,
\] 
showing that $G$ is quasicooperative.  After reindexing variables we
assume there is a fundamental subgraph $\Gam^1\subset \Gam (F)$ with
vertex set $\{1,\dots, n_1\}, \ 1\le n_1\le n$.  
Now apply Proposition \ref{th:fund0}.
\end{proof}

%%%%%%%%%%%%%%%%%%%%%%%%%%%%%%%%%%%%%%%%%%%%%%%%%%%%%%%%%%%%
\section{Monotone dynamics}   \mylabel{sec:mondyn}
%%%%%%%%%%%%%%%%%%%%%%%%%%%%%%%%%%%%%%%%%%%%%%%%%%%%%%%%%%%%
A local semiflow  $\Phi$ 
is {\em monotone} if $x\succeq y \implies \Phi_t x \succeq
\Phi_t y$.  Throughout this section we assume:
\begin{itemize}
\item {\em
 $\Phi:=\{\Phi_t\}_{t\ge 0}$ is a monotone local semiflow in  an
 ordered space $X$} 
 \end{itemize}
To simplify notation we may write $x(t):=\Phi_t x$ whenever $\Phi_t x$ is
defined. 
It is well known for the data in Equation (\ref{eq:basic}) that if $F$
is cooperative and $X$ is convex, the corresponding local semiflow
$\Phi$ is monotone.  This is a corollary of the M\"uller-Kamke theorem
\cite{Muller26, Kamke32} on differential inequalities (Hirsch
\cite{Hirsch82a}).

%%%%%%%%%%%%%%%%%%%%%%%%%%%%%%%%%%%%%%%%%%%%%%%%%%%%%%%%%%%%
\begin{proposition}		\mylabel{th:hs}
%%%%%%%%%%%%%%%%%%%%%%%%%%%%%%%%%%%%%%%%%%%%%%%%%%%%%%%%%%%%
The following are true for all $x\in X$:
\begin{description}

\item[(a)] No points of $\om x$ are related  by  $\rhd_X$ or $\lhd_X$

\item[(b)] $\om x$ is a singleton in the following cases:
\begin{description}

\item[(i)]  $\ov {\gam(x)}$ is compact and  there  exist $t_*\ge 0, \
  \varepsilon >0$ 
  such that 
\[ t_*<t<t_*+\varepsilon \implies 
  \Phi_t x \prec  x \text{ or } 
\Phi_t x \succ   x
\] 

\item[(ii)] $\ov {\gam(x)}$ is compact and there  exist $t> 0$
  such that  
\[
 \Phi_{t} x \rhd_X x \text{ or } \Phi_{t} x \lhd_X  x
\]
\end{description}
\end{description}
\end{proposition}
 %%%%%%%%%%%%%%%%%%%%%%%%%%%%%%%%%%%%%%%%%%%%%%%%%%%%%%%%%%%%
\begin{proof}
(a) and (b)(i) are sharpenings of Hirsch \& Smith \cite [Theorems
1.8, 1.4] {HirschSmith05}, respectively.   Assertion (b)(ii) follows
  from (b)(i). 
\end{proof} 
%%

%%%%%%%%%%%%%%%%%%%%%%%%%%%%%%%%%%%%%%%%%%%%%%%%%%%%%%%%%%%%
\begin{proposition}		\mylabel{th:gas0}
%%%%%%%%%%%%%%%%%%%%%%%%%%%%%%%%%%%%%%%%%%%%%%%%%%%%%%%%%%%%
Assume  $A\subset X$ is attracting.
\begin{description}

\item[(a)]  If each point  of $A$ is strongly accessible in
  $X$  from either above or below, then $A$ contains an equilibrium.
  
\item[(b)] If each point of $A$ is strongly accessible in $X$ from
  both above and below and $A\cap\E=p$ then $A=p$. 

\end{description}
\end{proposition}
%%%%%%%%%%%%%%%%%%%%%%%%%%%%%%%%%%%%%%%%%%%%%%%%%%%%%%%%%%%%
\begin{proof}  This is a slight generalization of  
  Hirsch  \cite[Theorems  III.3.1 and III.3.3] {Hirsch84},
  and the same proofs work here.  
 \end{proof}

%% {\tt *** From ord2.tex }

%%%%%%%%%%%%%%%%%%%%%%%%%%%%%%%%%%%%%%%%%%%%%%%%%%%%%%%%%%%%
 \begin{proposition}		\mylabel{th:ord2a}
%%%%%%%%%%%%%%%%%%%%%%%%%%%%%%%%%%%%%%%%%%%%%%%%%%%%%%%%%%%%
Assume $A\subset \om x$.  Let $q\in A$ be a minimal (respectively,
maximal) point of $A$ having a neighborhood $N\subset X$
such that there is a point $y\prec N$ (respectively, $y\succ N $) is
attracted to $A$.  Then $q=\inf A$ (respectively, $q=\sup A$).
\end{proposition}
%%%%%%%%%%%%%%%%%%%%%%%%%%%%%%%%%%%%%%%%%%%%%%%%%%%%%%%%%%%
\begin{proof}
To fix ideas we assume $q$ is a minimal point of  $A$ and $y \prec
N$.  
Notation is simplified  by  setting 
$\Phi_t w=w(t)$ whenevever $w\in X, t\ge 0$.   

Some point on $\gam(x)$ lies in $N$ its omega limit set contains $A$.  Replacing $x$ by such a point we
assume $x\in N$.  Therefore  $y\prec x$ and 
%%%%%%%%%%%%%%%%%%%%%%%%%%%%%%%%%%%%%%%%%%%%%%%%%%%%%%%%%%%%
\begin{equation}		\label{eq:tyx}
%%%%%%%%%%%%%%%%%%%%%%%%%%%%%%%%%%%%%%%%%%%%%%%%%%%%%%%%%%%%
y(t) \prec x(t),\qquad (t\ge 0)
\end{equation}
There is  a sequence  $t_n\to\infty$ such that 
 $x(t_n) \in N$ and  
%%%%%%%%%%%%%%%%%%%%%%%%%%%%%%%%%%%%%%%%%%%%%%%%%%%%%%%%%%%%
\begin{equation}		\label{eq:tyx2}
%%%%%%%%%%%%%%%%%%%%%%%%%%%%%%%%%%%%%%%%%%%%%%%%%%%%%%%%%%%%
x(t_n) \to q
\end{equation}
Because $\om y$ meets $A$ we can choose this sequence so that  also
%%%%%%%%%%%%%%%%%%%%%%%%%%%%%%%%%%%%%%%%%%%%%%%%%%%%%%%%%%%%
\begin{equation}		\label{eq:tyx3}
%%%%%%%%%%%%%%%%%%%%%%%%%%%%%%%%%%%%%%%%%%%%%%%%%%%%%%%%%%%%
  y(t_n) \to a\in A
\end{equation}
It follows from (\ref{eq:tyx}), (\ref{eq:tyx2}),
(\ref{eq:tyx3}) and closedness of the order relation that $a\preceq
q$,  so minimality of
$q$ implies $a=q$.  Thus
%%
%%%%%%%%%%%%%%%%%%%%%%%%%%%%%%%%%%%%%%%%%%%%%%%%%%%%%%%%%%%%
\begin{equation}		\label{eq:tyx4}
%%%%%%%%%%%%%%%%%%%%%%%%%%%%%%%%%%%%%%%%%%%%%%%%%%%%%%%%%%
 y(t_n)\to q \in \E
\end{equation}

Choose $n_0$ so that $y(t_{n_0}) \in N$.  If
$I\subset\Rp$ is a sufficiently small open interval about $t_{n_0}$ then
$s\in I\implies y(s)\in N$, hence $y(s)\succ y$.  The dual of Proposition
\ref{th:hs}(b)(i) now shows that $\om y$ is an equilibrium, hence $\om
y = \{q\}$.  It follows from (\ref{eq:tyx}) that $\om x \succ q$,
hence $ A \succeq q$.
\end{proof}

In the rest of this section we assume:

\begin{itemize}
\item 
{\em $X\subset\R n$ with the vector ordering. }
\end{itemize}
%%%%%%%%%%%%%%%%%%%%%%%%%%%%%%%%%%%%%%%%%%%%%%%%%%

%%%%%%%%%%%%%%%%%%%%%%%%%%%%%%%%%%%%%%%%%%%%%%%%%%%%%%%%%%%%
\begin{proposition}		\mylabel{th:moncon}
%%%%%%%%%%%%%%%%%%%%%%%%%%%%%%%%%%%%%%%%%%%%%%%%%%%%%%%%%%%%
Assume $x\in \om A, \ \om x = A$.  If 
$\inf A=p$ or $\sup A=p$ then $A=p$.  
\end{proposition}
%%%%%%%%%%%%%%%%%%%%%%%%%%%%%%%%%%%%%%%%%%%%%%%%%%%%%%%%%%%%
This result also holds when $X$ is ordered by a solid
polyhedral cone, but it is has not been proved for more general
ordered spaces.    For strongly order-preserving local semiflows a stronger
conclusion holds:  Every omega limit set is unordered (Hirsch \& Smith
\cite[Corollary 1.9] {HirschSmith05}).  
%% KEEP THE FOLLOWING SPACE %% 

\smallskip

\begin{proof} 
For any $\Sig\subset \{1,\dots,n\}$ the corresponding {\em face} of
 $\Rp^n$ is
\[
 \J:=\J(\Sig)=\{z\in \Rp^n\co z_i >0  \implies  i\in \Sig\}
\]
%%DA: little typo above,  z_i > 0
When $\Sig\ne\varnothing$ the corresponding {\em open face} is 
\[
  \Jo:=\Jo(\Sig)=\{z\in \Rp^n\co z_i >0 \dimply i\in \Sig\}
\]
It can be seen that $\ov {\Jo} =\J$  and $\Jo$ is
relatively open in its linear span.   Moreover
%%%%%%%%%%%%%%%%%%%%%%%%%%%%%%%%%%%%%%%%%%%%%%%%%%%%%%%%%%%%
\begin{equation}		\label{eq:j}
%%%%%%%%%%%%%%%%%%%%%%%%%%%%%%%%%%%%%%%%%%%%%%%%%%%%%%%%%%%%
(\forall z\in \Jo)\:(\exists \del >0) \quad  z\succ \J\cap
  N_\del (0)  
\end{equation}

Fix $x,p\in X$ such that 
$\inf \om x=p$ or 
$\sup \om x=p$; we have to prove $\om x=p$.  To fix ideas we assume
$p= 0 =\inf\om x$. 
{\em Claim: }
$\Phi_t(x)$ is defined for all $t \ge 0$.  It is well known that this
is the case if the orbit closure of $x$ is compact.  If it is not
compact, the orbit intersects the boundary of some open  ball
centered at $0$ in an infinite set.  Consequently $\om x$ contains a
point $\ne p$,  which implies the claim.     

For any $I\subset [0,\infty)$ set
$\Phi (I, x):=\{\Phi_t x\co t\in I\}$.  By the Baire category theorem
there is a dense open subset $S\subset [0,\infty)$ such that for each
component $I$ of $S$ there is a unique open face $\Jo_I\supset
\Phi (I, x)$.

There is a sequence $\{I_k\}$ of these components and points $t_k\in
I_k$ such that as $k\to \infty$ we  have
\[
   t_k \to \infty,\quad x(t_k) \succ 0, \quad x(t_k)\to 0
\]
%%DA: little typo: changed > to \succ above
  After  passing to a
subsequence we can assume there is an open face $\Ko$ such that
 $\Jo_{I_k}=\Ko$ for all $k$.  Choose such a $\Ko$ having the largest
possible dimension.  Then  $x(t)\in\Ko$ for
 sufficiently large $t$.  For if  $x(t_0) \in \Ko$ and  $\varepsilon
 >0$ is such that $x(t)
 \notin\Ko$ for $t\in (t_0, t_0 + \varepsilon]$, then $x(t')$ for some
 $t' \in(t_0, t_0 + \varepsilon]$ belongs to an open face of larger
 dimension, and this can only happen finitely many times. 
Set  $\dim \Ko = m\in \{1,\dots,n\}$ and relabel variables so that  
$\Ko = \Ko (\{1,\dots,m\})$. 

By (\ref{eq:j}) there exists $t_* >0$ such that 
\[
t >t_* \implies   x(t_*)\succ x(t) \succ 0 
\]  
By Proposition \ref{th:hs}(b)(i) the trajectory of $x(t_*)$
converges, necessarily to $0$.  Therefore $\om x=\om {x(t_*)}=
0$.
\end{proof}
%%

%%%%%%%%%%%%%%%%%%%%%%%%%%%%%%%%%%%%%%%%%%%%%%%%%%%%%%%%%%%%
\begin{corollary}		\mylabel{th:moncor}
%%%%%%%%%%%%%%%%%%%%%%%%%%%%%%%%%%%%%%%%%%%%%%%%%%%%%%%%%%%%
Assume $0\in X\subset\Rp^n$.   If $0\in \om x$, then $0=\om x$. \qed
\end{corollary}
%%%%%%%%%%%%%%%%%%%%%%%%%%%%%%%%%%%%%%%%%%%%%%%%%%%%%%%%%%%%
\begin{proof} Follows from Theorem \ref{th:moncon} because 
$0=\inf\om x$.
\end{proof}

\paragraph{Remark }
We digress to interpret this result biologically.  Let $x_i\ge 0$
stand for the ``size'' of species $i$ (population, biomass, density,
\dots) and call $\sum_{i=1}^n x_i$ the ``total size''.  Assume that
from each initial state $x(0)\in \Rp^n$ the species develop along a
curve $x(t)=(x_1(t),\dots, x_n(t))\in \Rp^n, t\ge 0$ governed by a
cooperative system (suggesting symbiosis or commensalism) in $\Rp^n$.
Then:
\begin{itemize}
\item  {\em
If the total population does not die out, the total size is bounded
above $0$.}  
\end{itemize}
This follows from the contrapositive of the Corollary.

\smallskip
The next result will be used to start the inductive proof of Theorem
 \ref{th:A}.  It applies only to cooperative systems, but the
 assumptions on $\Phi$, $X$ and $A$ are weaker than in Theorem
 \ref{th:A}.
Recall that every nonempty compact set in an ordered space contains  a
maximum point and a minimum point (Ward \cite {Ward54}),

%%%%%%%%%%%%%%%%%%%%%%%%%%%%%%%%%%%%%%%%%%%%%%%%%%%%%%%%%%%%
\begin{theorem}		\mylabel{th:fintrans}
%%%%%%%%%%%%%%%%%%%%%%%%%%%%%%%%%%%%%%%%%%%%%%%%%%%%%%%%%%%%
Assune $X\subset \R n$ has the vector ordering and $\Phi$ is a
monotone local semiflow in $X$.  Let  $A\subset X$ be  attracting and
finitely transitive for $\Phi$.   If every point of $A$ is strongly
accessible in $X$
from above or below, then $A\in\E$. 

\smallskip

More precisely: 
If $q\in A$ is maximal and strongly
accessible in $X$ from above then $A=q$.  Likewise if $q\in A$ is
minimal and strongly accessible in $X$ from below.
\end{theorem}
%%%%%%%%%%%%%%%%%%%%%%%%%%%%%%%%%%%%%%%%%%%%%%%%%%%%%%%%%%%%
\begin{proof}
It suffices to assume $q\in A$ is maximal and strongly accessible in
$X$ from above.  Under the current assumptions
there exist $x\in A, y\in X$ and neighborhood $U\subset X$ of $q$ such
that that $ q\in \om x, \ y\succ q$ and $y$ is attracted to $A$.
Evidently $q$ is maximal in $\om x$, hence $q=\sup \om x$ by
Proposition \ref{th:ord2a}, and therefore $q=\om{x}$ by Theorem
\ref{th:moncon}.

Suppose $z\in X$ and $\om z\cap U \ne\varnothing$.  There exists $l\in
 \Np$ with $z(t_l)\in U$, hence $y(t) \rhd_X z(t+l), \ (t\ge 0)$,
 and monotonicity proves
%%%%%%%%%%%%%%%%%%%%%%%%%%%%%%%%%%%%%%%%%%%%%%%%%%%%%%%%%%%%
\begin{equation}		\label{eq:dom}
%%%%%%%%%%%%%%%%%%%%%%%%%%%%%%%%%%%%%%%%%%%%%%%%%%%%%%%%%%%%
\om z\cap U \ne\varnothing \implies \om z \preceq q
\end{equation}
Now we prove for all $v\in X$:
%%%%%%%%%%%%%%%%%%%%%%%%%%%%%%%%%%%%%%%%%%%%%%%%%%%%%%%%%%%%
\begin{equation}		\label{eq:dom2}
%%%%%%%%%%%%%%%%%%%%%%%%%%%%%%%%%%%%%%%%%%%%%%%%%%%%%%%%%%%%
 q\in \om v \implies q=\om v
\end{equation}
For there exists $z\in \gam (v)\cap U$ and  Equation
(\ref{eq:dom}) implies $q=\sup \om z$,   hence $q=\om z=\om v$ by 
Theorem \ref{th:moncon}.

Let $\{a_k\}$ be any sequence in $U\cap A$ converging to $q$.  By
hypothesis there is a finite set $S\subset X$ such that each $a_k$ is
an omega limit point of some member of $S$.  By finiteness of $S$
there is a subsequence $\{b_k\}$ of $\{a_k\}$ and $v\in S$ such that
$\{b_k\}\subset \om v$.  Evidently $q\in \om v$, whence $q=\om v$ by
(\ref{eq:dom2}).  This can only happen if $b_k=q$ for all $k\in \Np$.
It follows that $q$ is isolated in the connected set $A$, entailing
$A=q$.
\end{proof}
\section{Proofs of the main theorems }   \mylabel{sec:proofs}
%%%%%%%%%%%%%%%%%%%%%%%%%%%%%%%%%%%%%%%%%%%%%%%%%%%%%%%%%%%%

%%%%%%%%%%%%%%%%%%%%%%%%%%%%%%%%%%%%%%%%%%%%%%%%%%%%%%%%%%%%
\paragraph{Proof of Theorem \ref{th:A}}
%%%%%%%%%%%%%%%%%%%%%%%%%%%%%%%%%%%%%%%%%%%%%%%%%%%%%%%%%%%%
Let the system $(F, X, \R n)$ be as in Theorem \ref{th:A}, with  a
finitely transitive attracting set $A\subset X$.  

{\em Step (i) } Consider first the case that $F$ is cooperative. Then
$\Phi$ is monotone because $X$ is convex, and each of the  assumptions (i),
(ii) implies each point of $A$ is strongly accessible in $X$ from
above or below. The conclusion for this case follows from Proposition
\ref{th:fintrans}.

\smallskip

{\em Step (ii) } We proceed by induction on $n$, the case $n=1$
following from the cooperative case.  Assume inductively that $n >1$
and that the conclusion holds for smaller values of $n$.  By Step (i) we
can assume $F$ is not cooperative, whence by Theorem \ref{th:cascade}
there is a cooperative system $(F^1, X^1, \R {n_1})$ and a cascade
$\Pi \co F \twoheadrightarrow F^1$ with $1\le n_1 <n$, whose fibre
systems are quasicooperative.  Lemma \ref{th:fibsys} shows that $(F_p,
X_p, E_p)$ is a fibre system for each $p\in\E (F^1)$.

The set $\Pi (A)\subset X^1$ is finitely transitive for the
cooperative system $F^1$, hence $\Pi (A)=p\in\E (F^1)$ by Step (i).
Thus $A$ lies in the invariant set $X_p= X\cap Q^{-1}(A)$, and $A$ is
attracting and finitely transitive for $\Phi_p:=\Phi|X_p$.  The
inductive hypothesis applied to $(F_p, X_p, E_p)$ shows that $A$ is an
equilibrium, completing the induction.  \qed

%%%%%%%%%%%%%%%%%%%%%%%%%%%%%%%%%%%%%%%%
\paragraph{Proof of Theorem \ref{th:B} }
%%%%%%%%%%%%%%%%%%%%%%%%%%%%%%%%%%%%%%%%
Consider first the case that $F$ is cooperative.
Assume {\em per contra} that the orbit closure of $x\in X$ contains a
nonempty open subset $U\subset X$.  As some open subset of $\R n$ is
dense in $X$ we can assume $U$ is open in $\R n$.  The orbit
$\gam(x)$, being a smooth curve,  is nowhere dense in $U$ because $n \ge 2$.
Therefore $U\subset \om x$, hence $\om x$ contains points $a, b$ such
that $a\rhd_X  b$.  But this contradicts Proposition \ref{th:hs}(a).

Now assume $F$ is not cooperative. By Theorem \ref{th:cascade} there
is a cascade $\Pi \co F \twoheadrightarrow F^1$ with $F^1$
cooperative.  If $W\subset X$ is open and $\gam$ is an orbit of $F$,
then $\Pi(W)$ is open in $X^1$ and $\Pi (\gam)$ is an orbit of $F^1$.
The cooperative case shows that $\Pi (\gam)\cap \Pi(W)$ is not dense
in $\Pi(W)$ 
and therefore $\gam\cap W$ is not dense in $W$.  \qed

%%%%%%%%%%%%%%%%%%%%%%%%%%%%%%%%%%%%%%%%%
\paragraph {Proof of Theorem \ref{th:C} }
%%%%%%%%%%%%%%%%%%%%%%%%%%%%%%%%%%%%%%%%%
If $F$ is cooperative, as when $n=1$, the conclusion follows from
Proposition \ref{th:gas0}. 
We proceed by induction on $n$, assuming that $n>1$ and the
theorem holds for smaller values.   

We can assume $F$ is not cooperative.  By Theorem \ref{th:cascade}
there is a cascade $\Pi \co F \twoheadrightarrow F^1$ cooperative
system $(F^1, X^1, \R {n_1})$ with $F^1$ cooperative and $1\le n_1
<n$, such that if $p\in \E(F^1)$ then $(F_p, X_p, E_p)$ is a
quasicooperative system.  Applying the inductive hypothesis twice, we
conclude that there exists $p\in \E(F^1)$ and $q\in \E (F_p)\subset
\E(F)$.

Assume $\E (F)= q$  and set $\Pi
(q)=p'\in\E (F^1)$.  Then $p'=p$.  For we showed above that every
fibre system contains an equilibrium of $F$, which must be $q$.  Thus
$\Pi^{-1} (p')$ and $\Pi^{-1}(p)$ are not disjoint, hence
they coincide and $\Pi$ maps both of them to $p$. 

By the inductive hypothesis $p$ is the global attractor for $F^1$,
therefore $X_p$ attracts all points of $X$ by Equation
(\ref{eq:piphi}).  This implies $A$ is the global attractor for
$\Phi|X_p$, and the inductive hypothesis applied to $(F_p, X_p, E_p)$
shows that $A=q$.  \qed
%%

%%%%%%%%%%%%%%%%%%%%%%%%%%%%%%%%%%%%%%%%%%%%%%%%%%%%%%%%%%%%
\section{Appendix:  Notes on the development of the concept  ``attractor''}
\mylabel{sec:attractors}
%%%%%%%%%%%%%%%%%%%%%%%%%%%%%%%%%%%%%%%%%%%%%%%%%%%%%%%%%%%%
\begin{quote}{\small \em
 In spite of the fact
 that everyone who is interested in dynamics has a more or less vague
 intuition of what an attractor of a map $f\colon M\rightarrow M$
 should be, there is no generally accepted mathematical definition for
 this concept even if $M$ is a smooth manifold and $f$ is also
 smooth.}
\hspace*{\fill}
 ---{\small H. Bothe \cite{Bothe87}}
\end{quote}
The first mathematical use of the word ``attractor'' may be in
Coddington \& Levinson's 1955 book \cite{CL55}, where it refers to an
asymptotically stable equilibrium.  The term was subsequently extended
to include an attracting cycles.  Today there are many definitions,
usually meaning an invariant set (of some kind) that is  approached uniformly
(in some sense) by the forward orbits of all (or most)
points in some neighborhood of the set.
 
Attractors do not occur explicitly
in the work of Poincar\'e or Birkhoff.  These authors were primarily
interested in Hamiltonian systems, which have no attractors because
they preserve volume.  

An early proof of existence of a unique attracting periodic orbit for
a general class of systems is in the 1942 paper of N. Levinson and
O. Smith \cite{LevinsonSmith42}.
 \footnote{
Thanks to George Sell for this reference.}

Early computer simulations revealed what appear to be attractors.  As
far back as 1952, Turing \cite{Tu52} published pictures of numerical
simulations of a nonlinear dynamical model of cell development,
exhibiting striking pattern formation.  Simulations by Stein \& Ulam
\cite{SU59, SU64} and Lorenz \cite{Lo63} gave persuasive pictorial
evidence of complicated structure in attractors, but attracted little
attention when they were published.  Hamming's review \cite{Ham65} of
\cite{SU64} was unenthusiastic:
\begin{quote} 
\small Many photographs of cathode ray tube displays are given, a
fondness for citing large numbers of iterations and machine time used
is revealed, and a crude classification of the limited results is
offered, but there appears to be no firm new results of general
mathematical interest\ldots

One can only wonder what will happen to mathematics if we allow the
undigested outputs of computers to fill our literature.  The present
paper shows only slight traces of any digestion of the computer
output.  
\end{quote}
\normalsize

Much of the early theoretical work on attractors on global analysis
was concerned with characterizing them in terms of Liapunov functions
and topological dynamics (e.g., Ura \cite {Ura53}, Auslander {\em et
al.}\ \cite {AusBhatSeib64}, Mendelson \cite {Me60}, Bhatia
\cite{Bh67}).  Little was known of their internal dynamics beyond the
existence of fixed points in global attractors for flows in Euclidean
space (Bhatia \& Szeg\"o \cite{BZ67}).

In the 1960s a number of articles on attractors and related forms of
stability were inspired by Sell \cite{Sell66}.  In his seminal 1967
work on global analysis, Smale gave detailed constructions
and analyses of hyperbolic attractors and other invariant sets, which
would later be called ``chaotic'' and ``fractal'', and  proved
them structurally stable.  He called attention to the vast mixture of
periodic, almost periodic, homoclinic and other phenomena found in
structurally stable attractors, even in rather simply given systems. 

``Strange attractors'' were proposed in 1971 as a model of turbulence
by Ruelle and Takens \cite{RuelleTakens71, RuelleTakens71a, Ruelle72},
Newhouse {\em et al.\ }\cite{NRT78}).  The physical significance of
this route to chaos is still  debated.

In his controversial 1972 book on morphogenesis (\cite{Thom72,
Thom75}) the late Ren\'e Thom issued a bold manifesto proclaiming the
fundamental scientific role of attractors: 
\begin{quote}
\small{ 
1.
Every object, or every physical form, can be represented by an {\em
attractor} $C$ of a dynamical system in a space $M$ of {\em internal
variables}.

\smallskip
\noindent
2. Such an object possesses no stability, and for this reason cannot be
perceived, unless the corresponding attractor is {\em structurally
stable}.

\smallskip
\noindent
3. Every creation or destruction of forms, every morphogenesis, can be
described by the disappearance of the attractors representing the
initial forms and their replacement through capture by the attractors
representing the final forms.  This process, called {\em catastrophe},
can be described in a space of {\em external variables}. \ldots}
\smallskip
\end{quote}

In recent years much work has been devoted to analysis of attractors
in specific classes of chaotic systems, such as those named after
Duffing, Lorenz, H\'enon and Chua, and to attractors having particular
topological properties, such as R. Williams' expanding attractors
(Williams \cite{Wi69}, Plykin \& Zhirov \cite{PZ93}).  A novel
measure-theoretic type of attractor due to Milnor \cite{Milnor85} has
stimulated several papers.

Many authors have investigated attractors in infinite-dimensional
systems, especially for partial differential equations, a prime
desideratum being finite dimensional gllobal attractors.  The large
literature includes books by Constantin {\em et al.\
}\cite{ConstFoiasNicol89}, Hale \cite {Hale88}, Ladyzhenskaya
\cite{Ladyzhen91}, Ruelle \cite{Ruelle95}, Sell \& You \cite{SY}, and
others.

Attractors, being objects defined  by topological limiting
processes,   resist classification and even description.  A general
theory appears quite distant.

\end{document}